\documentclass[article,noinfoline]{imsart}

\usepackage{amsthm,amsmath,amssymb,natbib}

\usepackage{graphicx}
\usepackage{color}

\startlocaldefs
\def\j{{\bf j}}
\def\k{{\bf k}}

\newcommand{\bdf}{\boldsymbol}
\endlocaldefs

\begin{document}

\begin{frontmatter}

% "Title of the paper"
\title{ Deconvolution, convex optimization, non-parametric empirical Bayes and treatment  
of non-response.  }
\runtitle{Deconvolution non-response}

% indicate corresponding author with \corref{}
% \author{\fnms{John} \snm{Smith}\corref{}\ead[label=e1]{smith@foo.com}\thanksref{t1}}
% \thankstext{t1}{Thanks to somebody}
% \address{line 1\\ line 2\\ printead{e1}}
% \affiliation{Some University}

%\author{\fnms{Lawrence D.} \snm{Brown,}\thanksref{t1}\ead[label=e1]{lbrown@wharton.upenn.edu}}
%\address{University of Pennsylvania; \printead{e1}}
%\thankstext{t1}{Research supported by an NSF grant.}
%\affiliation{University of Pennsylvania}
\author{\fnms{Eitan} \snm{Greenshtein}\ead[label=e2]{eitan.greenshtein@gmail.com}}
\address{Israel Central Bureau of Statistics; \printead{e2}}
\affiliation{Israel Central Bureau of Statistics}
\and
\author{\fnms{Theodor } \snm{Itskov} \ead[label=e3]{itsmatis@gmail.com}}
\address{   Israel Central Bureau of Statistics ; \printead{e3}}
%\thankstext{t3}{Research supported by an ISF grant.}
\affiliation{Israel Central Bureau of Statistics; \printead{e3}}

\runauthor{Greenshtein, Itskov}

\end{frontmatter}

{\bf Abstract.}  

Let $(Y_i,\theta_i)$, $i=1,...,n$, be independent random vectors distributed like $(Y,\theta) \sim G^*$, where the marginal distribution of $\theta$ is completely unknown, and
the conditional distribution of $Y$ conditional on $\theta$ is known. It is desired to estimate the marginal distribution of $\theta$ under $G^*$,
as well as functionals of the form $E_{G^*} h(Y,\theta)$ for a given $h$, based on
the observed $Y_1,...,Y_n$.

In this paper we suggest a deconvolution method for the above estimation problems
and discuss some of its
applications in Empirical Bayes analysis. The method involves a quadratic programming step, which is 
an elaboration on
the formulation and technique  in Efron(2013).  It is computationally efficient and may
handle large data sets, where the popular 
method, of deconvolution using EM-algorithm, is impractical. 

The main application that we study is  treatment of non-response.
Our approach  is nonstandard and does not involve  missing at random type of  assumptions.
The method is demonstrated in simulations, as well as in an analysis of a real data set from the
Labor force survey in Israel. Other applications including estimation of the risk,
and estimation of False Discovery Rates, are also discussed.

We also present a method, that involves convex optimization, for constructing confidence intervals for $E_{G^*} h$,
under the above setup.

\section{Introduction, Preliminaries and Examples.}

Consider a general empirical  Bayes setup, where
$(Y_i,\theta_i) $, are i.i.d., $i=1,...,n$, distributed like $(Y,\theta) \sim G^*$, 
 and the conditional distribution of $Y$ conditional on $\theta$ is $F_\theta$, $\theta \in \Theta$.    The marginal 
distribution of $\theta$
under $G^*$ is denoted  $G$. 
Suppose we only observe $Y_1,...,Y_n$, and
we should estimate the parameters
$\theta_1,...,\theta_n$.
It is often the case that $G$ is unknown and
should be estimated in the process of estimating the unknown parameters. 
We concentrate on the non-parametric empirical Bayes setup where
$G$ is completely unknown, 
as opposed to the parametric setup where $G$ is assumed to
be a member of some parametric family.

We have two main novel contributions in this paper.
One  is  suggesting a new deconvolution method, for the purpose of estimating $G$
by a corresponding  estimator $\hat{G}$. The deconvolution method
is based on quadratic
programming. Note, an estimator $\hat{G}$ for $G$ induces a corresponding 
estimator $\hat{G}^*$ for $G^*$, through  $d\hat{G}^*(y,{s})= dF_s(y)d\hat{G}(s) \equiv dG^*(y|\theta=s)d\hat{G}(s)$.
The other main contribution  is a nonstandard application of deconvolution and empirical Bayes
to the problem of treating non-response. Other applications are also described.

In the canonical examples of empirical Bayes the ultimate goal is the estimation of the parameters $\theta_i$, $i=1,...,n$,
based on the observed $Y_1,...,Y_n$.  However, our main interest and emphasis is on estimating
various functionals of the form $E_G h $ and $E_{G^*} h$ for various functions $h$.
We also consider the more general setup where  $(X_i,Y_i,\theta_i)$, $i=1,...,n$ are i.i.d., distributed like $(X,Y,\theta) \sim G^*$; the joint distribution, $G$, of $X$ and $\theta$ is completely unknown, while  $G^*(Y|X,\theta)$
the conditional distribution of $Y$ conditional of $X$ and $\theta$ is known.  We observe $n$ independent realizations
${\cal T}(X_i,Y_i)$, $i=1,...,n$ for some function ${\cal T}$ .  Here the pair $(X_i,\theta_i)$ may be
considered as  the `parameter' that determines the conditional distribution of $Y$, however, unlike  the former setup, the
"X-part" of the  `parameter'
is observed through  ${\cal T}(Y_i,X_i)$,
i.e., the parameter $(X,\theta)$ is not completely latent. The goal is again to estimate
$E_{G^*} h(X,Y,\theta)$ for various functions $h$. The estimators are of the form $E_{\hat{G}^*} h $, for a "deconvolution-estimator",
$d\hat{G}^*(x,y,s)=dG^*(y|\theta=s,X=x)d\hat{G}(x,s)$.

In Section 3 we present a method  for constructing a confidence interval for quantities  of the form $E_{G^*} h(X,Y,\theta)$,
based on ${\cal T} (X_i,Y_i) \; i=1,...,n$.
The main idea of that method is defining  an appropriate convex optimization problem, where the target function is linear and the
constraints are convex.

\bigskip

In the rest of this section we elaborate on a few empirical Bayes examples, where it is desired to estimate quantities
of the form $E_{G^*} h$. Our primary example is the problem of treating non-response.  Finally, we explain why
estimation of $G$ by $\hat{G}$ and then estimation of $E_G h$ by $E_{\hat{G}}h $ is a good alternative
to, say, mle estimation  of $\frac{1}{n} \sum h(\theta_i)$
by $\frac{1}{n} \sum h(\hat{\theta}_i)$, where $\hat{\theta}_i=\hat{\theta}_i(Y_i)$ is the
point-wise mle estimator of $\theta_i$, $i=1,...,n.$

\bigskip

\subsection{ Examples.}

{\it Deconvolution, Empirical-Bayes, and Estimation of the risk.}

In the canonical examples of Empirical Bayes, the ultimate goal is to estimate the individual parameters $\theta_i$, $i=1,...,n$.
In such problems the estimation of $E_{G^*} h$  for various $h$ could still  be central, as demonstrated in the following.

Let $\delta(Y)$ be a  decision function and $L(\theta,\delta(Y))$ a loss function.
Of a primary interest is the  quantity
\begin{equation} \label{eqn:risk}
 E_{G^*}L(\theta,\delta(Y))= E_G R(\theta,\delta)= E_G h_\delta(\theta);
\end{equation}
 here $Y \sim F_\theta$,  and $R(\theta,\delta) \equiv h_\delta(\theta)$ is the risk of $\delta$
conditional on $\theta$.
Thus, the quantity in $(\ref{eqn:risk})$ is the Bayes risk that corresponds to the decision function $\delta$,
under the loss $L$  and the prior $G$. 
The Bayes procedure is  $$\delta^B=argmin_\delta    E_G h_\delta(\theta) . $$
Uniformly good estimation of $ E_G h_\delta(\theta)$ over all $\delta$,
yields  good estimates of $\delta^B$. 

Once an estimator $\hat{G}$ for $G$ is obtained,
a natural approach is to let
\begin{equation}  \label{eqn:appl} \hat{\delta}^B=argmin_\delta    E_{\hat{G}} h_\delta(\theta) . \end{equation} 
Under a squared loss, the estimated decision function
in ( \ref{eqn:appl}) is $$ \hat{\delta}^B(y)= E_{\hat{G}^*}( \theta| Y=y).$$

More generally, under squared loss, in the setup where 
$(X_i,Y_i,\theta_i) \sim G^*$ and we observe ${\cal T}(X_i,Y_i)$,
$i=1,...,n$, a natural estimator for $\theta$ based on ${\cal T}(X,Y)$ is:
$$ \hat{\delta}( {\cal T}(X,Y))= E_{\hat{G}^*} (\theta| {\cal T} (X,Y)). $$
The case where $X$ and $Y$ are independent conditional on $\theta$ is of a special interest, e.g., as in our simulation section.

We should remark that ${\cal T}$ may be a randomized  transformation. In Brown et.al.
(2013), the set up is $(X_i,Y_i,\theta_i) \sim G^*$, $i=1,...,n$, are i.i.d, where the conditional distribution of $Y$ conditional on $\theta$ is $Poisson(\theta)$, while $X \sim Poisson (h)$ is independent of $\theta$ and $Y$. In that paper
even though a direct observation of $Y$ is available, the approach is to estimate the optimal decision function with respect to the artificially
"`corrupted"' observation
${\cal T}(X,Y)=X+Y$, with $h=h_n \rightarrow 0$, as $n \rightarrow \infty$. 
This approach is shown to have advantages relative to, say,
the classical EB estimator for a Poisson parameter, suggested by Robbins. In the sequel we will
not consider randomized  ${\cal T}$, although it is within our formulation.

There are common examples,
e.g., Poisson, Normal, where an estimator for $\delta^B$ may be obtained
directly without the estimation of ${G}$ and application of (\ref{eqn:appl}). On the direct approach for the estimation of $\delta^B$, versus approaches that
involve the estimation of $G$, see, e.g., Efron (2013). On direct
estimation of $\delta^B$ in the normal case see, e.g., Brown and Greenshtein (2009); on direct estimation in the Poisson case see, e.g., Brown, et.al. (2013). 

\bigskip

{\it Deconvolution and  Variations on False Discovery Rate}

Problems that involve estimation of $E_{G^*} h$ are related also to the problem of estimating false discovery rates (FDR), see, Benjamini and Hochberg (1995).

Let  $(Y_i,\theta_i) \sim G^*$, $i=1,...,n$, be independent  where conditional on $\theta_i$,
$Y_i \sim F_{\theta_i}$, $i=1,...,n$. Consider first the problem where it is desired to estimate
the proportion of indices $i$ for which $\theta_i> C$. When $n$ is large and $G$,
the marginal of $\theta$  is known, a reasonable trivial estimator, is $P_G(\theta>C)$.
Note that, $P_G(\theta>C) = E_G h(\theta)$
for the function $h$ which is the indicator of the event $\{ \theta>C \}$. When $G$ is unknown and estimated by $\hat{G}$,
the induced estimator is $E_{\hat{G}} h$. 

We now treat the more general (FDR) problem.
In order to fix ideas  consider the case $F_{\theta_i}=N(\theta_i,1)$. Suppose that  we consider observations  $i$ with $Y_i>A$ for some $A$ as "suspected discoveries", while
we consider as "true discoveries", observations for which $\theta_i >C$. In order to estimate the proportion of true
discoveries among suspected discoveries we should estimate the quantity $E_{G^*} h$  for the function $h$ which is
the indicator of the event $((\theta_i>C )\cap ( Y_i>A))$. 

When it is desired to estimate the proportion of "true discoveries" among suspected discoveries for a given realization, the following perspective and 
alternative approach might be beneficial.  Let $G^{*t}$ be the conditional distribution of $(Y,\theta)$ conditional on the event $Y>A$, we  treat the observations
$(Y_i,\theta_i)$ with $Y_i \leq A$, as truncated and the remaining ones are treated as  i.i.d., observations from $G^{*t}$, where $G^{*t}(y|\theta)=F_\theta(y|Y>A)$.
Let $G^t$ be the marginal distribution of $\theta$ under $G^{*t}$ and $\hat {G}^t$ its "deconvolution-estimate", let $h$ be an indicator of the event $\theta>C$,
we may estimate  $E_{G^t} h$ by $E_{\hat{G}^t} h$. 
See Greenshtein et. al. (2008),  for treatment of a related  problem.

\bigskip

{\it Deconvolution  and Treatment of Non-Response.}
A  main novel contribution of this paper, is
an application of our deconvolution method to treat non-response. The proposed treatment of non-response
does not involve the, often assumed and seldom verifiable, assumption
of Missing At Random (MAR), conditional on some covariates.

Let ${\bf S}=\{ i_1,...,i_n \}$ be a random set of indices 
that correspond to randomly sampled items from a finite population
of size $N$, indexed by $\{1,...,N\}$. Those are the indices of the items 
in the population who i) were randomly sampled for a survey ii)  responded.

Suppose, it is desired to estimate the total 
$T= \sum_{i=1}^N X_i$ in the population, based on the $n$
available observations. Let $I_i$ be an indicator of the event
"`item $i \in {\bf S}$"', $i=1,...,N$. 
Let $p_i = E(I_i)$. 
Then $$\hat{T}=\sum_{i \in {\bf S} }\frac{X_i}{p_i},$$
is the Horvitz Thompson estimator for $T$. It is an unbiased
estimator, as may be seen from the representation
\begin{equation}
 \label{eqn:0-1}
 \hat{T}= \sum_{i=1}^N \frac{X_i}{p_i}I_i.
 \end{equation}  

Typically, $p_i$, $i=1,...,N$ are unknown although the sampling probabilities are known. This is since the 
corresponding response probabilities are unknown. Thus, the above estimator can not be applied.

We will approximate (\ref{eqn:0-1}), by a nonparametric empirical
Bayes modeling together with a deconvolution step.
Consider a situation where there is an additional covariate
$Y_i$ for every item $i$, $i \in {\bf S}$, such that $Y_i \sim F_{p_i}$.
In one example, that we will give, $Y_i$ is the number
of visits until a response was obtained;  
in another example
$Y_i$
is the number of responses of item $i$
in a longitudinal panel survey, where  each sampled  item is
attempted to be interviewed four times in four consecutive months. 
  
We model the observations $(X_i,Y_i,p_i)$, $i \in \bdf{S}$ as i.i.d   $ (X_i,Y_i,p_i) \sim G^{*t}$. 
Here $G^{*t}$  is the conditional distribution of $(X_i,Y_i,p_i)$ conditional on $ i \in \bf{S}$, which is different than
$G^*$ the distribution of $(X_i,Y_i,p_i)$, $i=1,...,N$.
A  natural estimator for (\ref{eqn:0-1}) is  $n E_{G^{*t}} \frac{X}{p} \equiv nE_{G^{*t}} h(X,p)$, for
$h(X,p)=\frac{X}{p}$. This treatment is under truncation, where we have no knowledge about the observations that correspond to indices $i$,
$i \notin \bf{S}$ . Under censoring, when there exists partial information
about items with index $i$, $i \notin  {\bf S}$, related ideas will be applied. The formal distinction and different treatment under truncation versus censoring
will be explained and demonstrated in Sections 4,5.

A general reference to sampling is, e.g., Lohr(2009). A reference for missing data and non-response issues is, e.g., Little and Rubin (2002).

\bigskip

{\it Non parametric maximum likelihood estimation of $G$}

The first study of the estimation of $G$, under the above setup,
 was conducted by Kiefer and Wolfowitz (1956). They suggested to find the 
non-parametric mle
for $G$, and also gave conditions under which the non-parametric mle estimator
$\hat{G}$ converges weakly to the true $G$. Estimation of $ E_G h$
by $E_{\hat{G}} h$ is often much better than estimating the individual
parameters $\theta_i$, say by an mle, $\hat{\theta}_i(Y_i)$, and then average, to obtain the estimator   $\sum^n h(\hat{\theta}_i)/n$.
This is demonstrated in the following Example 1.
\bigskip

{\bf Example 1.} Consider the Normal example where 
$F_{\theta_i} = N(\theta_i,1)$. 
Let $h$ be the  function $h(\theta)= 1/\theta$.
Suppose it is known that the support
of $G$ is bounded bellow by 0.5, but otherwise it is completely unknown.  Then the mle for $\theta_i$ is $\hat{\theta}_i=max(0.5,Y_i)$. Now suppose that the true $G$ has a point mass at 1. By Kiefer and Wolfowitz (1956), the mle
$\hat{G}$ for $G$ converges weakly to $G$, so 
$E_{\hat{G}} h \rightarrow 1$. However, a quick simulation 
shows that $\frac{1}{n}\sum \frac{1}{\hat{\theta}_i}$ converges to
1.19.  

\bigskip

Estimation of sums of the form $\sum_i h(Y_i,\theta_i)$, was studied by
Zhang (2005). Further examples may be found there, as well as a study of the efficiency of certain estimators.

\bigskip

The rest of the paper is organized as follows.
In Section 2, our deconvolution method is explained. In Section 3, we present a method that involves convex optimization,
to construct confidence intervals for quantities of the form $E_{G^*}h$.
In Section 4 we present `empirical Bayes type Horvitz Thompson' estimators
in the context of treating non-response. In Section 5, the derivation and performance of those estimators are illustrated through  a simulated
practical example. In Section 6, we demonstrate
our method for treating non-response, through a real data set from the Labor Force Survey in Israel.

\section{Deconvolution using quadratic programming.}
\label{subsec:cond}

In this section we present a deconvolution algorithm which 
involves  quadratic programming.

Our deconvolution is a method for deriving a Non-Parametric  Maximum Likelihood Estimator (NPMLE) for a `prior' $G$.
We use the term deconvolution in a wide sense, that includes identifying mixtures, as studied, e.g., by Lindsay (1995), Lindsay and Roeder (1993), Lee et.al., (2013), and, of course, the fore mentioned seminal paper of Kiefer and Wolfowitz (1956).
Our quadratic programming approach, rather than the more common EM-algorithm, is in line with the general suggestion and advocation
of Koenker and Mizera (2013) for the usage of convex optimization.
It may be applied on high dimensional problems with tens of thousands of
observations and general mixing $G$, where the complexity of EM
algorithms makes them impractical.

In the following subsection we treat the problem of estimating the marginal distribution of a latent-variable/unknown-parameter.  The same ideas apply  to the more
general problem of  estimating the joint distribution of a latent variable and an observed variable. We present the ideas in two stages where the general case is formulated in 
subsection (\ref{subsec:uncond}).

Our approach is based  on the setup and formulation in
Efron (2013). We elaborate more by defining and solving an appropriate quadratic programming problem.

\subsection{ Deconvolution for the estimation of the marginal distribution of a latent variable.} \label{subsec:uncond}

Consider a standard empirical Bayes setup,
as described in the introduction, where $(Y_i,\theta_i) \sim G^*$, are i.i.d., $i=1,...,n$. We
 assume  discrete distributions,  in particular $F_{\theta}, \; \theta \in \Theta$, are discrete
with a common finite support denoted $\{ y_1,...,y_J\}$,   and $G$ is discrete with 
a given support $\{s_1,...,s_K \}$. The treatment of the continuous
cases may be done through discretization. In principle the discretization of the $Y$-variables should be more delicate as the number of observations 
increases, but the `right' way of discretization is beyond the scope of this paper. Our main examples and applications 
in sections 4-6 involve discrete observations $Y_i$, $i=1,...,n$, specifically, censored Geometric and Binomial.
The considerations that are involved in the discretezation of $G$ have to do with the complexity of the estimation algorithm.

Our observations $Y_i$, $i=1,...,n$, are independent
and identically 
distributed like a random variable $Y$. Denote their discrete density by ${\bdf f}=(f_1,...,f_J)'$, where
$f_j=P(Y=y_j), \; j=1,...,J$. Denote
$$ p_{jk}= P(Y=y_j|\theta=s_k), j=1,...,J; \; k=1,...,K.$$  Denote the density of the discrete distribution $G$ by ${\bdf g}=(g_1,...,g_K)'$, where
$g_k=P_G(\theta=s_k), \; k=1,...,K$.
Denote by $P$  the ${J \times K}$ matrix $P=(p_{jk})$.

Then:  \begin{equation} \bdf{f}= P \bdf{g}. \end{equation}
 
Recall, the support of $G$ is known (or practically approximated
by a dense grid $\{s_1,...,s_K\}$ ), it is the density $\bdf{g}$ that
should be estimated. 
We now reduce the problem through sufficiency.  Note that a sufficient
statistic is $\hat{\bdf{f}}=(\hat{f}_1,...,\hat{f}_J)'$,  where 
$\hat{f}_j$ is the proportion of observations among $Y_1,...,Y_n$,
that had the value $y_j$, $j=1,...,J$.  
Now, $\bdf{\hat{f}}$ is a scaled multinomial vector with 
mean $\bdf{f}$ and a corresponding covariance matrix $\Sigma_{\bdf{f}}/n$. Its distribution is asymptotically
multivariate normal.  Note, that there is a linear dependence,
thus the corresponding covariance matrix  $\Sigma_{\bdf{f}}^{-1}$ does not exist. We may  replace 
$\bdf{\hat{f}}$ by  the sufficient statistic
$\bdf{\hat{f}^*}=(\hat{f}_1,...,\hat{f}_{J-1})'$,
whose corresponding covariance matrix is  $\Sigma^*/n$.
The mean of $\bdf{\hat{f}^*}$ is $P^*\bdf{g}$, where $P^*_{(J-1)\times K}$ is obtained from $P$ by deleting its last column.
Since the distribution of $\bdf{\hat{f}^*}$ is asymptotically multivariate normal, a  solution $\hat{\bdf{g}}$ to:

\begin{equation}
 \label{eqn:mle} \min_{ \bdf{g}} ( \bdf{ \hat{f}^{*}}-P^*\bdf{g} )' \Sigma^{*-1} ( \bdf{\hat{f}^{*}}-P^*\bdf{g}), \end{equation}

$\;\;\;\;\;\;$  s.t. $0 \leq g_i \leq 1$, $\sum g_i=1$,

is asymptotically an mle estimator for $\bdf{g}$.
Note!, we write `an mle' rather than `the mle' since a
solution and an mle are not necessarily unique. See, also 
Remark 1 bellow.  Practically, $\Sigma^*$
is replaced  by its estimate, which is obtained by utilizing 
the  multinomial distribution of $n\hat{\bdf{f}}$. A  special care should be 
taken when estimating  $\Sigma^{*-1}$, since $\Sigma^*$ might
be close to being singular. 
Our approach
in our numerical work was to add $0.001$ to the diagonal of the regular estimator of the covariance matrix of a multinomial vector, then take its inverse as the estimator of $\Sigma^{*-1}$.

\bigskip

{\it  Calibration.}

Suppose there is a function $A$, for which it is known that
$E_G A(\theta)=a$. In such a case we may add 
to the above linear programming the linear
constraint:

$$ \sum g_k A(s_k)=a.$$
Similarly, when there are a few such functions $A_1,...,A_b$.

\bigskip

\bigskip

The numeric work in this paper was done by applying the 
quadratic programming 
function {\it ipop}, from the R-package {\it kernlab}, Karatzoglou, et. al. (2004).

\bigskip

{\bf Remark 1:} 

It may be concluded from
Lindsay and Roeder (1993) or Lindsay (1995) that when there are
only $J$ possible values to $Y$, there exists an  mle for $G$, that  has  $J-1$ points or less
in its support. Thus,  we can not expect consistency of 
an arbitrary mle estimator $\hat{G}$, unless the support of  $G$ is known to have no more
than $J-1$ points. However, if $Y$ is obtained by a discretization
of a continuous observation,
which may become more and more delicate as $n$ grows, we may expect
consistency when $J=J_n \rightarrow \infty$. 

Furthermore, by adding calibration constraints, we might
get an mle which has a larger support and the corresponding estimator $\hat{G}$ 
is a
better approximation
of $G$. 
\bigskip

\bigskip

\subsection{ Deconvolution for estimation of the joint distribution of a
latent and an observed variables.}
\label{subsec:uncond}

In the previous section we considered the problem of estimating the distribution
$G$ of a latent variable $\theta$. In this section we will generalize the method to estimate the joint distribution of a latent variable $\theta$ and an
observed variable $X$, where $(X,Y,\theta) \sim G^*$. Let 
$(X_i,Y_i, \theta_i) \sim G^*$, $i=1,...,n$ be independent. We only observe  
${\cal T}(X_i,Y_i)$, $i=1,...,n$ for some ${\cal T}$, and the estimation is based only on those observed values. 

The variables $X_i$ are discrete, their possible values are $x_1,x_2,...,x_L$.

%Consider the following  model for the joint distribution $G^*$ 
%of $(X,Y,\theta)$. The marginal probability of the event $X=x_l$ is $\alpha_l$, $l=1,...,L$. Conditional upon $X=x_l$, $\theta$ is distributed $G_l$, 
%$l=1,...,L$. Finally conditional upon $X$ and $\theta$, $Y$ is distributed $F_\theta$. The distributions  $G_l$ 
%and the probabilities $\alpha_l$, $l=1,...,L$ are completely unknown.
%We observe $n$ independent  pairs $(X_i,Y_i)$ realized by the above mechanism. 

Our goal is to estimate the joint distribution of $\theta$ and $X$, which is determined by

$$ g_{lk}=P_{G^*}(X=x_l,\theta=s_k), \; l=1,...,L, \; k=1,...,K,$$ 
we denote 
$ {\bf g}=(g_{11}, g_{12},...,g_{LK}) \equiv (g_1,g_2,...,g_{L \times K}),$
note the dual indexing of the vector $\bf{g}$.

Let $t_1,...,t_Q$, be the distinct  values of ${\cal T}(x_l,y_j), \; l=1,...,L, \; j=1,...,J$. 
We assume that the conditional distribution $G^{*}({\cal T}=t| X=x, \theta=s)$ is known, thus the
$L\times K$ pairs ${\bf v}_{lk}=(x_l, s_k)$, $l=1,...,L, \; k=1,...,K$, play the role of the "`parameter"' that governs the conditional distribution.
Denote by ${\bf v}=( {\bf v}_{11}, {\bf v}_{1,2},...,{\bf v}_{LK}) \equiv ({\bf v}_1,...,{\bf v}_{L \times K})$, the vector of "parameters", note the dual
indexing of the vector ${\bf v}$.

As  in the previous subsection,   denote by
$$ p_{\j \k} = P({\cal  T}=t_\j| (X,\theta)={\bf v}_\k) , \; \j=1,...,Q, \; \k=1,..., L \times K,$$
let $P=(p_{\j\k})$ be the corresponding matrix as in the previous subsection.

Given $n$ observations, let ${f}_\j=P({\cal T}(X_i,Y_i)=t_\j)$, $\j=1,...,Q$,  let ${\bf f}=(f_1,...,f_Q)'$,
then $${\bf f}=P {\bf g}.$$
Let $\hat{f}_{\j}$, $\j=1,...,Q$ be the proportion of observations $i$ for which ${\cal T}(X_i,Y_i)=t_\j$, let
 $\hat{\bf{f}}^*=(\hat{f}_1,...,\hat{f}_{Q-1})$. Then ${\bf f}^*  \equiv E\hat{ {\bf f}}^* = P^* {\bf g}$, for the matrix $P^*$, which is obtained from $P$, as in the previous subsection,
by deleting its last row.

Let   $\Sigma^*/n$ be the covariance matrix of $\hat{\bdf{f}}^*$, and suppose it is non-singular. Then,

\begin{equation}
\label{eqn:mle1} \min_{ \bdf{g}} ( \bdf{ \hat{f}^{*} }-P^*\bdf{g} )' \Sigma^{*-1} ( \bdf{\hat{f}^{*}}-P^*\bdf{g}), \end{equation}

$\;\;\;\;\;\;$  s.t. $0 \leq g_{lk} \leq 1$, $\sum_{lk} g_{lk}=1,$

is asymptotically an mle estimator for $\bdf{g}$.

\bigskip

%The estimator  $\bdf{\hat{g}}$ induces an estimator $\bdf{\hat{g}_l}$ for the 
%conditional densities $\bdf{g}_l=(g_l^1,...,g_l^k)'$,  $l=1,...,L$,
%$g_l^i =g_{li}/\alpha_l$,  by the obvious
%normalization letting
%\beg%in{equation} \label{eqn:cnd}
%\hat{g}_l^i= \frac{\hat{g}_{li}}{\sum_j  \hat{g}_{lj} }, \; i=1,...,k.
%\end{equation}

%The above defines the following estimator for 
%$E_{G_l} h \equiv E_{G^*}(h | X=x_l)$,
%the conditional expectation of $h$ under $G_l$, 
%\begin{equation}
%\label{eqn:condex}
%E_{\hat{G}_l} h = \sum_{i=1}^k \hat{g}_l^i h(s_i), \; i=1,...,k.
%\end{equation}  

\bigskip

{\it Calibration}

The above quadratic programming may incorporate various additional linear constraints. For example suppose that there is an indicator $I=I(X)$, $I=1$
if the corresponding measurement was taken from a male, $I=0$
otherwise. Suppose it is known that $P_{G^*}(I=1)=0.5$. Then the constraint
$$\sum_{(l,k):I(X_l)=1} g_{lk}=0.5,$$ may be added to the quadratic programming
defined in  (\ref{eqn:mle1}).

\section{ Confidence intervals and linear optimization.}
\label{sec:conf}

We consider the setup of the previous section, where we observe i.i.d
$(X_i,Y_i, \theta_i) \sim G^*$, $i=1,...,n$.

Suppose it is desired to estimate the expectation 
$$T=E_{G^*} h(X,\theta) =\sum_{k,l} (x_l,s_k) g_{lk}.$$
Note, a simple modification of the treatment  bellow  applies also for expectations of the form
$T=E_{G^*} h(X,Y,\theta)$; however in order to simplify the notations we consider the above functionals.

It is of an interest to obtain a confidence interval for $T$, this could reassure that an mle estimator (recall, often the mle is  not unique), is giving
a reliable estimate.

Let $\hat{\bf{f}}^*$  and $\Sigma^*$ be  as in the previous section.
Suppose that $\Sigma^*$ is non-singular. Let $\hat{\Sigma}^*$ be the empirical 
covariance matrix. Then as the sample size approaches infinity 
$\hat{\Sigma}^{*-1}$
approaches $\Sigma^{*-1}$ in probability. Furthermore, 
the distribution of $\sqrt{n} \hat{\bf{f}}^*$
converges weakly  to a multivariate normal distribution with covariance matrix
$\Sigma^*$. Recall, under the general setup of subsection 3.2, we observe ${\cal T} (X,Y)$, whose support is of
size $Q$.

Consider the solution of  the following problem, of linear optimization under convex
constraints.
\begin{eqnarray} \label{eqn:LP}
&& T_U=\max_g  \sum_{l,k} h(x_l,s_k) g_{lk} \\
&& T_L=\min_g  \sum_{l,k} h(x_l,s_k) g_{lk} \nonumber
\end{eqnarray}

s.t.
 
$$ n ( \bdf{ \hat{f}^* }-P^*\bdf{g} )' \hat{\Sigma}^{*-1} ( \bdf{\hat{f}^{*}}-P^*\bdf{g}) < \chi^2_{(Q-1),1-\alpha},$$

$$\;\;\;\;\;\;   0 \leq g_{lk} \leq 1, \;\;  \sum_{l,k} g_{lk}=1,$$
in the above $\chi^2_{(Q-1),1-\alpha}$ is the critical value of the appropriate $\alpha$-level $\chi^2$ test with $Q-1$ degrees of freedom. As before, additional 
convex calibration constraints nay be added if available.

\bigskip

{\bf Theorem 1:} If $\Sigma^*$ is non-singular, then $(T_L,T_U)$ is 
 a $(1-\alpha)$
level confidence interval, asymptotically as $n \rightarrow \infty$.

\bigskip

The above theorem is for discrete  variables $\theta$, $X$ and $Y$. For
continuous cases a discretization should be done. The general guide lines for discretization is that $Q$  will be  of size $o(n)$, so that 
there will be enough observations in each of the $Q-1$ "`cells"' 
and the asymptotic $\chi^2_{Q-1}$ distribution will  hold; the considerations involved in the
discretization of $\theta$ and $X$ have to do with the complexity of the convex optimization. A formal asymptotic treatment of the discretization is beyond the scope of this paper.

\section{ Non-response and Empirical Bayes type Horvitz Thompson estimators.}

A general survey from a population of size $N$ indexed by $\{1,...,N\}$, may be described as follows. Each subject $i$,
$i=1,...,N$,
 in the population
is sampled with probability $\pi_i$ for an interview, but once subject $i$ is sampled  a response from that subject is obtained with
probability $p^*_i \leq 1$.  Let $\bdf{S}$ be the random set of indices, corresponding to subjects who  i) were sampled for an interview ii) responded.
Then, for subject $i$, $P( i \in \bdf{S})= \pi_i p^*_i=p_i$, $i=1,...,N$.
We define the indicator random variable $I_i$, $I_i=1$ iff $i \in  \bf{S}$; denote
$P(I_i=1)=p_i$.

In many surveys  the subjects are equally likely to be sampled to the survey, i.e., $\pi_i \equiv \pi$ are all equal. In the following we treat this case. Thus, w.l.o.g., we may assume that $\pi=1$,  and $p_i=p^*_i$. Modification of the treatment bellow applies when
$\pi_i$, $i=1,...,N$, may have `a few' possible values.

%Now, suppose that each item has a measurement of interest $X_i$, $i=1,...,N$, $X_i$ is discrete and may have the values
%$l, \; l=1,...,L$. Let $\bdf{S}_l$, $l=1,...,L$  be the random set of indices who i) selected ii) responded iii) their $X$-value is $l$.
%Let $T=\sum^N_{i=1} X_i$. Then, assuming that $0 < p_0 \leq min_i p_i$, the Horvitz Thompson estimator  for $T$ is:
%$$ \hat{T}=\sum_l  l \sum_{i \in \bdf{S}_l} \frac{1}{p_i}.$$
%The last estimator may be used by an  oracle, who knows the values of $p_i$.

We will apply our deconvolution technique and the
empirical Bayes ideas, to provide Horvitz Thompson type of estimators in the context of Empirical Bayes. 

We model the items of the size $N$ population, as realizations of $N$, i.i.d random vectors $(X_i,Y_i,p_i,I_i)$, which are distributed like 
$(X,Y,p,I) \sim G^*$, $X$ is the variable of interest. 
 The joint distribution of $X$ and $p$ is arbitrary
and the conditional distribution of $Y$  conditional on $X$ and $p$ is known;
\begin{equation} \label{eqn:condit}  p_i= P(I_i=1)=P(I_i=1|X_i,p_i). \end{equation}
In order to fix ideas think of $X_i$
as an employment-status of item $i$, $I_i$ indicator of the event "`item $i$
was sampled for a survey and responded"'. 
In one of our examples in the sequel, 
$Y_i$ is the number of attempts until a response was obtained from 
subject $i$,  where there are at most
$M_0$ attempts. Thus, in this example $I_i=0$ iff $Y_i>M_0$.
We model 
$Y \sim F_p=Geometric(\tilde{p})$ for $\tilde{p}=\tilde{p}(p)$,
$p=1-(1-\tilde{p})^{M_0}$.

\bigskip
{\it Truncated versus censored observations.} We will consider two different setups. In one setup the event $I_i=0$ means
that the observation is truncated, i.e., we do not know about variables
with $I_i=0$, and thus, our available observations may be considered as an i.i.d sample
from the distribution, denoted $G^{*t}$, of $(X,Y,p,I)$, conditional on $I=1$.
Another setup is of censored observations where we do know
about the event  $I_i=0$; e.g., in the example where $Y_i$ is the number of visits
until a response, the event $I_i=0$ implies $Y_i>M_0$. The two setups lead to    
two versions of our general deconvolution technique, in the truncated setup we estimate the
joint distribution of $p$ and $X$ under $G^{*t}$, while in the censored setup
we estimate the joint distribution of $p$ and $X$ under $G^*$.  The joint distribution of $X$ and $p$ under
$G^*$ and $G^{*t}$ will be denoted by $G$ and by $G^t$ correspondingly.

 \subsection{ Empirical Bayes type Horvitz Thompson estimators.}

Suppose we want to estimate $T= E\sum_{i=1}^N X_i$. We now present three unbiased estimators. Those are in fact pseudo-estimators since they are functions of the unknown $p_i$, however they will be modified  later to become legitimate estimators. 

$$ \hat{T}_0= \sum \frac{X_i}{p_i} I_i \equiv \sum X_i I_i A_0^i ,$$

$$\hat{T}_1=\sum X_i I_i E( \frac{1}{p_i} | X_i, I_i=1) \equiv \sum X_i I_i A_1^i ,$$

$$\hat{T}_2= \sum \frac{X_i}{E(p_i| X_i)} I_i \equiv   \sum X_i I_i A_2^i. $$

In the above $A_j^i$ are implicitly defined. The estimator $\hat{T}_0$ is
basically the standard Horvitz Thompson estimator.

\bigskip

{\bf Theorem 2:}

i) Under the condition $p_i>0$ w.p.1,  $E(\hat{T}_1)=E(\hat{T}_0)=T$. 
 Under the (weaker)
condition $E(p|X)>0$  w.p.1, $E(\hat{T}_2)=T$.

ii) Under the condition $p_i>0$ w.p.1, $Var(\hat{T}_1) \leq Var(\hat{T}_0)$.

iii) Under the condition $p_i>0$ w.p.1, $\hat{T}_2=\hat{T}_1$.

 \bigskip

{\bf Proof:} We prove the theorem for the case $N=1$.

i) $E\hat{T}_0=T$ follows immediately, similarly to the implication for a standard Horvitz-Thompson estimator.
$E\hat{T}_1= T$ follows since
$$ E E(\hat{T}_0| X,I) =E \hat{T}_1.$$

Assume that $E(p|X)>0$, w.p.1, then
$$E(\hat{T}_2) = E\frac{X I}{E(p|X)} =  EE( \frac{X I}{E(p|X)} | X)= EX=T,$$
the third equality in the above follows since by (\ref{eqn:condit}) $E(I|X)=E(p|X)$.

ii) The assertion follows by a Rao-Blackwell argument, due to the above
conditional expectation representation.

iii)  The assertion follows since   
$$ dG^*(p|I=1,X_i=x)=\frac{ p dG^*(p|X_i=x)}{ \int p dG^*(p|X=x_i)},$$
whence $A_1^i= \int \frac{1}{p} dG^*(p|I=1,X=x_i)= \frac{1}{E(p|X=x_i)}
=A_2^i$.

\bigskip

In practice the terms $A_j^i$, $j=1,2$ are unknown, they will be estimated using our
deconvolution method through the estimation of the joint distribution of
$X$ and $p$, under $G^*$ and $G^{*t}$ respectively, for $j=1,2$.

For every $i, \; i=1,...,n$, define 
$$\tilde{A}_2^i= 1/E_{\hat{G}} (p|X=X_i); $$ here $\hat{G}$ is the deconvolution estimator for $G$, the joint distribution of $X$  and $Y$ under $G^*$.
%$E_{\hat{G}_l} h $ is given in (\ref{eqn:condex} ). Here $\hat{G}_l$, 
%estimate the conditional distribution of $p$ conditional on $X=x_l$, 
%$l=1,...,L$, under
%$G^{*}$, i.e., the corresponding $\hat{g}_{li} $ estimate the joint probability of $X$ and $p$
%under $G^{*t}$. 

In the truncated setup the role of $G^*$ in our deconvolution method is played
by the conditional distribution $G^{*t}$.  Now, $$\tilde{A}_1^i= E_{\hat{G}^t}(\frac{1}{p}|X=X_i);$$
here $\hat{G}^t$ is the estimated joint distribution of $X$ and $p$ under $G^{*t}$.

%In order to emphasize it, we denote
%the estimated conditional distribution distribution of $p$ 
%conditional on $X=x_l$ under $G^{*t}$ by $G^t_l$. $l=1,...,L$.
%For every $i, \; i=1,...,n$, define $\tilde{A}_1^i= E_{\hat{G}_l^t}(\frac{1}{p})$when $X_i=x_l$; here $\hat{G}_l^t$ is  given in (\ref{eqn:cnd} ) and
%$E_{\hat{G}_l^t} h $ is given in (\ref{eqn:condex} ). 
%We re-emphasize, the distribution $G^{*t}$ plays the role of $G^*$, i.e., $\hat{g}_{li} \equiv \hat{g}^t_{li} $ estimate the joint probability of $X$ and $p$
%under $G^{*t}$. 

We now present the legitimate versions of $\hat{T}_1$ and $\hat{T}_2$, i.e., 
estimators which are functions only of the available observations,

\begin{equation} \label{eqn:tild1}
\tilde{T}_1=\sum X_i I_i \tilde{A}_1^i,
\end{equation}

\begin{equation} \label{eqn:tild2}
\tilde{T}_2=\sum X_i I_i \tilde{A}_2^i.
\end{equation}

We relate the estimators $\tilde{T}_1$ and $\tilde{T}_2$ through the Horvitz-Thompson estimator. However, in fact
the estimation of $T$ under the censored  setup may be done without the mediation of the Horvitz-Thompson estimator.
 In fact, an mle estimator for $T=N E(X)$  under the censored setup is: $$ \tilde{T}_3= N E_{\hat G} X= N\sum_l x_l \sum_k \hat{g}_{lk},$$ for a corresponding,
mle, estimator ${\bf {\hat g} }$.  Asymptotically $\tilde{T}_2 \approx \tilde{T}_3$. This  may be seen, by the following.
Denote by $m_l$, $l=1,...,L$, the number of indices $i$  satisfying $X_i=x_l$, then $E  m_l =  N E (p|X=x_l) \sum_k {g}_{lk}$, whence for large $N$,
$m_l \approx  N E_{\hat{G}} (p|X=x_l) \sum_k {\hat g}_{lk}$; note that
$\tilde{T}_2=\sum_l x_l m_l/E_{\hat{G}}(p|X=x_l)$.
The later version, $\tilde{T}_3$, is better suited
compared to $\tilde{T}_2$, for deriving a confidence interval for $T$ by the method that is given in Section \ref{sec:conf}.

There are a few advantages to $\tilde{T}_2$ compared to $\tilde{T}_1$, the obvious one
is that it is defined also when the event $p=0$ has a positive probability.
The other advantage is since that in the estimation of
$A_2$ we use some additional censored information, which is not available
(i.e., truncated), in the estimation of $A_1$. 
Avoiding possible near  singularity for small $p$, involved in the estimation of $A_1^i$,
is another advantage in attempting to estimate $A_2^i$ when possible.
Finally, typically there is an available external  information about the 
distribution $G^*$, that may be used through calibration, while that information is typically unknown under $G^{*t}$.
\bigskip

In the following simulation sections, we will apply our estimators  in the estimation of the expected proportion $\alpha _{x_l}$ of  items with a corresponding
$X=x_l$.
Their expected total number is estimated  by 
\begin{equation}
\tilde{T}_{x_l}^j = \sum_{i: X_i=x_l} \tilde{A}_j^i,
\end{equation}
for $j=1,2$, for the truncated and censored setups correspondingly. The following formula applies for the
truncated and censored  estimators for the proportion $\alpha_{x_l}$ 
when setting $j=1,2$ correspondingly,

\begin{equation} \label{eqn:htalpha}
\hat{\alpha}^j_{x_{l_0}}= \frac{  \sum_{i: X_i=x_{l_0}} \tilde{A}_j^i  }
{ \sum_l \sum_{i: X_i=x_l} \tilde{A}_j^i  }.
\end{equation}

In the next section we will use as a benchmark the following estimator,
that could be used by an `oracle' that knows $p_i \; i=1,...,N$. Such an oracle
could estimate the size of the population with corresponding $X=x_l$, by
$ \sum_{i: X_i=x_l} \frac{ I_i}{p_i}$. The corresponding oracle estimator
for $\alpha_{x_{l_0}}$ would be:
\begin{equation}
\label{eqn:orc}
\text{oracle}_{x_{l_0}}=\frac{  \sum_{i: X_i=x_{l_0}} \frac{ I_i}{p_i}  }
{ \sum_l\sum_{i: X_i=x_l} \frac{ I_i}{p_i}  }
\end{equation}

\section{Simulations}

Consider a survey where in its first stage an initial subset of the population is sampled and in the next stage 
there is an attempt to interview each sampled subject.
As mentioned, we assume that each subject in the population is equally likely to be sampled in the first stage, with sampling probability  $\pi_i \equiv \pi, \; i=1,...,N$; w.l.o.g $\pi=1$. 

Suppose our policy is to make at most $M_0$ attempts in order to obtain a response from a sampled subject, however obviously if a response is obtained in the $j<M_0$ attempt, no further attempts are made. We model the number of attempts
until a response is obtained by subject $i$,  by
a Geometric random variable  with a success probability $\tilde{p}_i$, 
$i=1,...,N$; assume 
$0<\min_i  \tilde{p}_i$.  Let $Y_i$ denote the number of attempts until a response was obtained.

Assuming $\pi_i \equiv 1$, the probability $p_i$ of subject $i$,
$i=1,...,N$
to be in the set $\bdf{S}$, of items that i) were sampled for the survey
and ii) responded, is:

\begin{equation}
\label{eqn:map} 
p_i = 1- (1-\tilde{p}_i)^{M_0}. 
\end{equation}

Thus, there is a one to one correspondence between  $\tilde{p}_i$ and 
$p_i=P(I_i=1)= P(I_i=1|X_i,p_i)=P(Y_i \leq M_0|X_i,p_i)$.

\subsection{Truncated setup}
We are interested in the estimation of the joint distribution
of $X$ and $p$ under $G^{*t}$, i.e., conditional upon $I=1$. 

Note, the distribution of $Y_i$, conditional on $i \in \bdf{S}$ is:
\begin{equation}
\label{eqn:cens}
P(Y=j|p)= \frac{ (1-\tilde{p})^{j-1}\tilde{p}  }
{ 1- (1-\tilde{p})^{M_0} }, \; j=1,...M_0;
\end{equation}
here $\tilde{p}=\tilde{p}(p)$, as given in  (\ref{eqn:map}).

Denote the distribution of $Y_i$, $i \in \bdf{S}$,
given in (\ref{eqn:cens}) by $F_{p_i}$.

Given a grid of points $\{ s_1,...,s_K\}$
we define the vector ${\bf v}=( (x_1,s_1),(x_1,s_2),...,(x_L,s_K))=({\bf v}_1,...,{\bf v}_{L \times K})$. The possible outcomes, 
${\cal T}(X,Y)=(X,Y)$ are denoted  $((x_1,1),(x_1,2),...,(x_L,{M_0}))=
(t_1,...,t_Q)$, $Q=L \times M_0$. As in subsection (\ref{subsec:uncond}), let
${\bf f}=(f_1,...,f_Q)$, where $f_{\j}=P({\cal T}(X,Y))=t_{\j}$, $\j=1,...,Q$. For every $\j$,
$\j=1,...,Q$, we denote $t_\j=(t_{\j1},t_{\j2})$, 
for every $\k$, $\k=1,...,K \times L$,  ${\bf v}_\k=(v_{\k1},v_{\k2}).$

Denote  
$$ p_{\j \k}=P({\cal T}=t_\j| (X,p)={\bf v}_\k), \; \j=1,...,Q, \; \k=1,...,L \times K.$$
This defines the matrix $P=(p_{\j\k})$  as explained in
the previous section,

\[ p_{\j\k}= \left \{ \begin{array}{ll}
0 &  v_{\k1} \neq t_{\j1} \\
P(Y=t_{\j2}| p=v_{\k2}) &  v_{\k1} = t_{\k1} .
\end{array}
\right. \]

Note, ${\bf f}= P{\bf g} $. We proceed as in subsection (\ref{subsec:uncond})
to derive an estimator for ${\bf g}$. In turn we obtain the estimators 
$\hat{\alpha}^1_{x_l}$, $l=1,...,L$, as in (\ref{eqn:htalpha}).

%Now, suppose for every item $i$, $i \in \bdf{S}$ there is a corresponding
%value $X_i$, where $X_i$ may equal 0 or 1.
%Let $P^*$ be the matrix that is obtained from $P$ as explained in
%subsection \ref{subsec:cond}.

\subsection{Censored setup.}

We will repeat the estimation of $\alpha_{x_l}$, $x_l=0,1$,  under the same setup,  estimating the proportions by $\hat{\alpha}_{x_l}^2$,
i.e, the censored  version of  (\ref{eqn:htalpha}) .

In the current setup we observe ${\cal T}(X_i,Y_i)$, where  

\[ {\cal T}(X_i,Y_i)= \left\{ \begin{array}{ll}
(X_i,Y_i)  &   Y_i \leq M_0 \\
NR  &     Y_i>M_0
\end{array}
\right. \] 

Here $"NR"$  abbreviate "Non-Response" and the outcome $NR$ implies that $Y_i>M_0$.
We denote the possible outcomes by $(t_1,...,t_Q)=( (X_1,1), (X_1,2),...,(X_L,{M_0}),"NR")$.  
 The number of possible values of ${\cal T}$ is $Q=(M_0 \times L) +1$.
As in subsection (\ref{subsec:uncond}), let $\bf{f}^*$ be the vector of expected proportions of  the $Q-1$ possible outcomes when excluding the outcome "NR";
let ${\bf v}$ be  the $L \times K$ dimensional vector, as in  
(\ref{subsec:uncond}).

%\equiv ({\bf u}_1, {\bf u}_2,..., {\bf u}_{M\times k})$, note the dual indexing.
%Let  $p_{ji}= P( {\cal T}= {\bf u}_j | (X,p)={\bf v}_i)$. 
We write ${\bf v}_\k \equiv (v_{\k1},v_{\k2} )$,    ${t}_\j \equiv (t_{\j1},t_{\j2})$.

\[ p_{\j \k}= \left \{ \begin{array}{ll}
0 &  v_{\k 1} \neq t_{\j 1} \\
(1-\tilde{p}_{\k})^{t_{\j 2}}\tilde{p}_{\k} &  v_{\k 1} = t_{\j 1} .
\end{array}
\right. \]

Here $\tilde{p}_{\k}=\tilde{p}_{\k}(v_{\k2})$, is the probability of success of the Geometric random variable $Y$, while $v_{\k2}$ is the probability of
success within $M_0$ trials, $v_{\k2}=1-(1-\tilde{p}_k)^{M_0}$, as explained in the previous subsection.

We proceed as  in subsection (\ref{subsec:uncond}), to obtain the deconvolution estimator  for ${\bf g}=( g_{11},g_{12},...,g_{LK})$ 
that determines $G$, the joint
distribution of $X$ and $p$ under $G^*$.  The estimated $\hat{G}$ defines the estimator $\hat{\alpha}^2_{x_l}$, $l=1,...,L$, for $\alpha_{x_l}$, as explained
in ( \ref{eqn:htalpha}).

\subsection{Numerical experiments}

\bigskip

In the following we simulate populations of size
$N$, where we randomly assigned to $N^0 \sim Binomial(N,0.5)$ items
a corresponding value $X=0$, and to the remaining $N^1$ a corresponding $X=1$ was assigned. A value $\tilde{p}$,
of a response probability in a single attempt, was 
randomly assigned to each of the $N^0$ items independently,
under a distribution $\tilde{G}_0$. Similarly a value $\tilde{p}$ was assigned randomly to each of the $N^1$ items based on a distribution
$\tilde{G}_1$. 
A corresponding pair $(G_0,G_1)$, of distributions  of
the possible values of response probabilities  $p$ is determined. 
Let $\alpha_{x_l}=0.5$, $x_l=0,1$, be the expected proportion of 
items with a corresponding $X=x_l$. Finally, for each item $i$, $i=1,...,N$,
a Geometric random  variable $Y_i \sim Geometric(\tilde{p}_i)$  was simulated.

We simulated scenarios with $N=1000$ and $N=10000$. 
The cases $M_0=4,6,8$, were studied for each of the
following three classes of pairs of distributions $(\tilde{G}_0,\tilde{G}_1)$,
parametrized by $\gamma$.

{\bf Two Points.}  The distribution $\tilde{G}_0$  has a two points support,
at the points $0.5 $ and $0.9$, with probability  mass 0.5 at each.

The distribution $\tilde{G}_1 \equiv \tilde{G}_1^\gamma $ is a $(-\gamma)$ translation of $\tilde{G}_0$.
We present results for the cases $\gamma= 0.1, 0.2, 0.3, 0.4$.

\bigskip

{\bf Uniform.} The distribution $\tilde{G}_0$ is uniform on the interval $(0.1,1)$.
The distribution $\tilde{G}_1 \equiv \tilde{G}_1^\gamma$, is a mixture of $\tilde{G}_0$ and a point mass at $0.1$, where the mixing weights are $(1-\gamma)$ and $\gamma$ correspondingly. We present results for $\gamma=0.1, 0.2, 0.3, 0.4$. 

\bigskip

{\bf Normal.} The distribution $\tilde{G}_0$ is a $N(0.5,0.1)$,  
`rounded up' to
$0.1$ and `rounded down' to $1$. The distribution $\tilde{G}_1 \equiv \tilde{G}_1^\gamma$ is 
$N(0.5-\gamma, 0.1)$ `rounded up' to $0.1$ and `rounded down' to $1$.
We present results for $\gamma=0.1, 0.2, 0.3, 0.4$. 

\bigskip

In the simulations we compared the performance of the following
estimators for $\alpha_{0}=0.5$.
The naive estimator, that   estimates $\alpha_{0}$ by  the sample proportion, i.e., 
the proportion among responders, of items $i$  with $X_i=0$; the
estimators  $\hat{\alpha}_{0}^j, \; j=1,2$, that correspond to the truncated and censored setups, 
as given  in (\ref{eqn:htalpha}) ; the 
`oracle' estimator as in (\ref{eqn:orc}). 

The grid points $\{s_1,...,s_K \}$ taken as the support of $G$, are induced by
$\{ \tilde{s}_1=0.1,\tilde{s}_2=0.12,...,\tilde{s}_K=1 \}$ that were taken as the support of $\tilde{G}$.

The following two tables correspond to the cases $N=1000$ and $N=10000$.
The columns S-naive, S-$\hat{\alpha}^1_0$, S-$\hat{\alpha}^2_0$,
S-oracle correspond to the square root of the simulated
mean squared error of each of the corresponding methods, based on 1000
repetitions.  The columns m-naive, m-$\hat{\alpha}^1_0$, m-$\hat{\alpha}^2_0$,
correspond to the simulated average of each of the corresponding methods; the simulated mean of the oracle's estimator was virtually $0.5$ and thus not presented.

It may be seen that $\hat{\alpha}^2_0$, clearly dominates $\hat{\alpha}^1_0$,
as may be expected. The performance of all the methods is improved by an increase in $M_0$, but the improvement is much sharper for $\hat{\alpha}^j_0, \; j=1,2$. For $\gamma=0$, the setup would become Missing at Random, in which
the naive estimator is the best. As $\gamma$ increases, the other methods
dominate the naive. It may be seen that for large enough $M_0$ and $\gamma$, in
all of our simulated configurations $\hat{\alpha}^2_0$, dominates the naive.
In the uniform case when $N=1000$, $\hat{\alpha}^1_0$ does not dominate the naive in any of the configurations, however when we let $M_0=10$, we get domination of 
$\hat{\alpha}^1_0$, specifically, for $N=1000$, $ M_0=10$ , $\gamma=0.4$,
S-$\hat{\alpha}^1_0$=0.0279, compared to S-naive=0.0383.

The performance of the estimator $\hat{\alpha}^2_0$, is comparable to that of
the oracle when $M_0=8$, and it is amazingly close to it in the Two-Points case.
It may be seen that our methods reduce the bias. This is important
beyond the reduction of the mse, since often the estimators 
arrive as a time-series and the final estimators involve additional smoothing of the time-series. Obviously, smoothing around the true value gives further 
reduction in mse, compared to smoothing of a  biased sequence.

Finally, an important  `moral' from  the two tables is that an increase in $M_0$
is much more important for risk reduction, relative to an increase in the sample size.  For example,
in the setup of Two-Points, $\gamma=0.4$, $M_0=6$, $N=1000$,  we have  
S-$\hat{\alpha}^1_0$=0.0274, while for $\gamma=0.4$, $M_0=4$, $N=10000$,
the mse is increased and
S-$\hat{\alpha}^1_0$=0.0478. Obviously, the number of interviewing attempts
in the first case is smaller than that in the second case.

  % Table generated by Excel2LaTeX from sheet 'Sheet1'

\begin{table}[htbp]
  \centering
  \caption{N=1000}
    \begin{tabular}{rrrrrrrrrr}
  %  \toprule
    {$\tilde{G}$} & {$M_0$} & {$\gamma$} & \text{m-naive} & 
		\text{m-$\hat\alpha_0^1$} & \text{m-$\hat\alpha_0^2$} & \text{S-naive} & \text{S-oracle} & \text{S-$\hat\alpha_0^1$} & \text{S-$\hat{\alpha}_0^2$} \\
%    \midrule
    TwoPts & 4     & 0.1   & 0.4909 & 0.4206 & 0.4963 & 0.0184 & 0.0161 & 0.0868 & 0.0186 \\
    TwoPts & 4     & 0.2   & 0.4743 & 0.4094 & 0.4891 & 0.0304 & 0.0165 & 0.0996 & 0.0250 \\
    TwoPts & 4     & 0.3   & 0.4470 & 0.3867 & 0.4766 & 0.0556 & 0.0173 & 0.1230 & 0.0405 \\
    TwoPts & 4     & 0.4   & 0.3978 & 0.3456 & 0.4532 & 0.1035 & 0.0192 & 0.1618 & 0.0668 \\
    TwoPts & 6     & 0.1   & 0.4966 & 0.4815 & 0.4995 & 0.0164 & 0.0161 & 0.0300 & 0.0164 \\
    TwoPts & 6     & 0.2   & 0.4872 & 0.4823 & 0.4978 & 0.0208 & 0.0165 & 0.0335 & 0.0173 \\
    TwoPts & 6     & 0.3   & 0.4663 & 0.4726 & 0.4937 & 0.0373 & 0.0162 & 0.0417 & 0.0203 \\
    TwoPts & 6     & 0.4   & 0.4221 & 0.4358 & 0.4846 & 0.0796 & 0.0178 & 0.0731 & 0.0274 \\
    TwoPts & 8     & 0.1   & 0.4978 & 0.4975 & 0.4992 & 0.0156 & 0.0154 & 0.0172 & 0.0155 \\
    TwoPts & 8     & 0.2   & 0.4933 & 0.5007 & 0.4996 & 0.0173 & 0.0160 & 0.0200 & 0.0160 \\
    TwoPts & 8     & 0.3   & 0.4788 & 0.5022 & 0.4990 & 0.0268 & 0.0165 & 0.0245 & 0.0169 \\
    TwoPts & 8     & 0.4   & 0.4394 & 0.4762 & 0.4948 & 0.0629 & 0.0178 & 0.0349 & 0.0185 \\
    Uniform & 4     & 0.1   & 0.4855 & 0.3739 & 0.4921 & 0.0224 & 0.0181 & 0.1335 & 0.0446 \\
    Uniform & 4     & 0.2   & 0.4682 & 0.3638 & 0.4816 & 0.0360 & 0.0184 & 0.1435 & 0.0548 \\
    Uniform & 4     & 0.3   & 0.4504 & 0.3562 & 0.4777 & 0.0530 & 0.0201 & 0.1516 & 0.0609 \\
    Uniform & 4     & 0.4   & 0.4301 & 0.3509 & 0.4710 & 0.0720 & 0.0197 & 0.1571 & 0.0664 \\
    Uniform & 6     & 0.1   & 0.4882 & 0.4441 & 0.4952 & 0.0205 & 0.0174 & 0.0629 & 0.0287 \\
    Uniform & 6     & 0.2   & 0.4738 & 0.4399 & 0.4893 & 0.0312 & 0.0174 & 0.0679 & 0.0340 \\
    Uniform & 6     & 0.3   & 0.4597 & 0.4347 & 0.4860 & 0.0438 & 0.0176 & 0.0735 & 0.0371 \\
    Uniform & 6     & 0.4   & 0.4457 & 0.4314 & 0.4858 & 0.0570 & 0.0183 & 0.0770 & 0.0388 \\
    Uniform & 8     & 0.1   & 0.4908 & 0.4757 & 0.4973 & 0.0189 & 0.0166 & 0.0340 & 0.0224 \\
    Uniform & 8     & 0.2   & 0.4794 & 0.4709 & 0.4941 & 0.0261 & 0.0162 & 0.0373 & 0.0238 \\
    Uniform & 8     & 0.3   & 0.4679 & 0.4687 & 0.4937 & 0.0362 & 0.0172 & 0.0408 & 0.0256 \\
    Uniform & 8     & 0.4   & 0.4555 & 0.4634 & 0.4913 & 0.0476 & 0.0173 & 0.0449 & 0.0255 \\
    Normal & 4     & 0.1   & 0.4792 & 0.3570 & 0.4966 & 0.0267 & 0.0168 & 0.1492 & 0.0227 \\
    Normal & 4     & 0.2   & 0.4422 & 0.3485 & 0.4917 & 0.0602 & 0.0176 & 0.1594 & 0.0295 \\
    Normal & 4     & 0.3   & 0.3859 & 0.3414 & 0.4863 & 0.1156 & 0.0199 & 0.1679 & 0.0404 \\
    Normal & 4     & 0.4   & 0.3231 & 0.3332 & 0.4833 & 0.1778 & 0.0211 & 0.1738 & 0.0471 \\
    Normal & 6     & 0.1   & 0.4902 & 0.4571 & 0.4989 & 0.0195 & 0.0169 & 0.0523 & 0.0184 \\
    Normal & 6     & 0.2   & 0.4664 & 0.4489 & 0.4955 & 0.0375 & 0.0169 & 0.0631 & 0.0214 \\
    Normal & 6     & 0.3   & 0.4223 & 0.4380 & 0.4920 & 0.0796 & 0.0180 & 0.0744 & 0.0257 \\
    Normal & 6     & 0.4   & 0.3691 & 0.4333 & 0.4919 & 0.1321 & 0.0191 & 0.0782 & 0.0272 \\
    Normal & 8     & 0.1   & 0.4945 & 0.4899 & 0.4987 & 0.0169 & 0.0160 & 0.0232 & 0.0168 \\
    Normal & 8     & 0.2   & 0.4777 & 0.4875 & 0.4968 & 0.0277 & 0.0166 & 0.0277 & 0.0178 \\
    Normal & 8     & 0.3   & 0.4461 & 0.4813 & 0.4964 & 0.0564 & 0.0170 & 0.0350 & 0.0187 \\
    Normal & 8     & 0.4   & 0.4016 & 0.4762 & 0.4962 & 0.0999 & 0.0183 & 0.0401 & 0.0196 \\
%    \bottomrule
    \end{tabular}%
  \label{tab:addlabel}%
\end{table}%

% Table generated by Excel2LaTeX from sheet 'Sheet1'
\begin{table}[htbp]
  \centering
  \caption{ $N=10000$}
    \begin{tabular}{rrrrrrrrrr}
%    \toprule
    {$\tilde{G}$} & {$M_0$} & {$\gamma$} & \text{m-naive} & {m-$\hat{\alpha}_0^1$} & {m-$\hat\alpha_0^2$} & \text{S-naive} & \text{S-oracle} & {S-$\hat\alpha_0^1$} &{S-$\hat\alpha_0^2$} \\
%    \midrule
    TwoPts & 4     & 0.1   & 0.4907 & 0.4191 & 0.4974 & 0.0106 & 0.0052 & 0.0837 & 0.0081 \\
    TwoPts & 4     & 0.2   & 0.4747 & 0.4119 & 0.4931 & 0.0258 & 0.0054 & 0.0925 & 0.0136 \\
    TwoPts & 4     & 0.3   & 0.4469 & 0.3939 & 0.4849 & 0.0534 & 0.0053 & 0.1110 & 0.0261 \\
    TwoPts & 4     & 0.4   & 0.3983 & 0.3478 & 0.4648 & 0.1019 & 0.0061 & 0.1559 & 0.0478 \\
    TwoPts & 6     & 0.1   & 0.4957 & 0.4786 & 0.4993 & 0.0067 & 0.0051 & 0.0237 & 0.0053 \\
    TwoPts & 6     & 0.2   & 0.4870 & 0.4786 & 0.4990 & 0.0139 & 0.0050 & 0.0261 & 0.0058 \\
    TwoPts & 6     & 0.3   & 0.4663 & 0.4767 & 0.4976 & 0.0341 & 0.0053 & 0.0303 & 0.0071 \\
    TwoPts & 6     & 0.4   & 0.4228 & 0.4459 & 0.4928 & 0.0774 & 0.0057 & 0.0574 & 0.0117 \\
    TwoPts & 8     & 0.1   & 0.4985 & 0.4979 & 0.5001 & 0.0052 & 0.0049 & 0.0066 & 0.0050 \\
    TwoPts & 8     & 0.2   & 0.4933 & 0.4976 & 0.4999 & 0.0084 & 0.0050 & 0.0088 & 0.0051 \\
    TwoPts & 8     & 0.3   & 0.4785 & 0.4999 & 0.4995 & 0.0221 & 0.0053 & 0.0126 & 0.0054 \\
    TwoPts & 8     & 0.4   & 0.4395 & 0.4829 & 0.4978 & 0.0607 & 0.0055 & 0.0206 & 0.0060 \\
    Uniform & 4     & 0.1   & 0.4845 & 0.3666 & 0.4926 & 0.0164 & 0.0057 & 0.1361 & 0.0321 \\
    Uniform & 4     & 0.2   & 0.4679 & 0.3617 & 0.4852 & 0.0326 & 0.0060 & 0.1413 & 0.0393 \\
    Uniform & 4     & 0.3   & 0.4504 & 0.3580 & 0.4833 & 0.0499 & 0.0062 & 0.1457 & 0.0437 \\
    Uniform & 4     & 0.4   & 0.4317 & 0.3524 & 0.4807 & 0.0686 & 0.0061 & 0.1516 & 0.0493 \\
    Uniform & 6     & 0.1   & 0.4874 & 0.4412 & 0.4945 & 0.0136 & 0.0052 & 0.0612 & 0.0193 \\
    Uniform & 6     & 0.2   & 0.4741 & 0.4363 & 0.4919 & 0.0264 & 0.0054 & 0.0663 & 0.0224 \\
    Uniform & 6     & 0.3   & 0.4600 & 0.4338 & 0.4910 & 0.0404 & 0.0057 & 0.0693 & 0.0238 \\
    Uniform & 6     & 0.4   & 0.4453 & 0.4299 & 0.4891 & 0.0550 & 0.0057 & 0.0736 & 0.0257 \\
    Uniform & 8     & 0.1   & 0.4897 & 0.4724 & 0.4976 & 0.0116 & 0.0054 & 0.0296 & 0.0119 \\
    Uniform & 8     & 0.2   & 0.4791 & 0.4692 & 0.4957 & 0.0216 & 0.0053 & 0.0331 & 0.0134 \\
    Uniform & 8     & 0.3   & 0.4681 & 0.4677 & 0.4958 & 0.0323 & 0.0054 & 0.0349 & 0.0140 \\
    Uniform & 8     & 0.4   & 0.4563 & 0.4653 & 0.4947 & 0.0440 & 0.0057 & 0.0375 & 0.0144 \\
    Normal & 4     & 0.1   & 0.4792 & 0.3498 & 0.4963 & 0.0215 & 0.0052 & 0.1519 & 0.0124 \\
    Normal & 4     & 0.2   & 0.4432 & 0.3398 & 0.4913 & 0.0571 & 0.0056 & 0.1633 & 0.0196 \\
    Normal & 4     & 0.3   & 0.3863 & 0.3287 & 0.4852 & 0.1138 & 0.0060 & 0.1749 & 0.0296 \\
    Normal & 4     & 0.4   & 0.3223 & 0.3276 & 0.4870 & 0.1777 & 0.0064 & 0.1751 & 0.0297 \\
    Normal & 6     & 0.1   & 0.4892 & 0.4553 & 0.4981 & 0.0120 & 0.0052 & 0.0466 & 0.0076 \\
    Normal & 6     & 0.2   & 0.4656 & 0.4456 & 0.4943 & 0.0348 & 0.0052 & 0.0574 & 0.0119 \\
    Normal & 6     & 0.3   & 0.4224 & 0.4330 & 0.4915 & 0.0778 & 0.0055 & 0.0710 & 0.0155 \\
    Normal & 6     & 0.4   & 0.3695 & 0.4275 & 0.4920 & 0.1306 & 0.0059 & 0.0756 & 0.0148 \\
    Normal & 8     & 0.1   & 0.4942 & 0.4897 & 0.4990 & 0.0077 & 0.0051 & 0.0131 & 0.0058 \\
    Normal & 8     & 0.2   & 0.4786 & 0.4851 & 0.4973 & 0.0221 & 0.0052 & 0.0186 & 0.0072 \\
    Normal & 8     & 0.3   & 0.4460 & 0.4791 & 0.4965 & 0.0543 & 0.0055 & 0.0260 & 0.0080 \\
    Normal & 8     & 0.4   & 0.4030 & 0.4743 & 0.4964 & 0.0971 & 0.0056 & 0.0300 & 0.0081 \\
%    \bottomrule
    \end{tabular}%
  \label{tab:addlabel}%
\end{table}%

\newpage

\section{ Analysis of real data of Labor Force Survey.}  
  
In this section we will apply our method on a real data set from the
Labor Force Survey, that is conducted by the Israel Central Bureau of Statistics.  The sampling method is 4-8-4 rotating panels, however
for our analysis, it may be equivalently treated and described as
a 4-in rotation, which is described in the following.

The survey is given to four panels, where each panel is investigated
for four consecutive months. Each month one panel finishes its fourth
investigation and in the next month it will be replaced by a new panel that will remain for four months. The main purpose of the survey is to estimate the proportion of  `Unemployment', `Employment',
and those who are `Not in Working Force (NWF)', the last category is of those who do not have a job nor they are looking for one.
Denote the corresponding values of our variable $X$-`working status', by
0, 1, 2.  We are interested in estimating $\alpha_0, \alpha_1$ and $\alpha_2$.
The population 
of interest is  of residents whose  age is above 15, and the proportions are with respect to that population. 
The probability $\pi$ to be included in the sample is the same for
each person. As explained, for our purpose of estimating
proportions we  assume w.l.o.g that $\pi=1$.

Temporarily assume that, we have only the data from the panel that is investigated for the fourth
time (`fourth panel'). Its size is about 5000, 
however, only $n$ responses were
obtained, $m_l$ responses from people with working status $x_l$, $x_l=0,1,2$. 
The general response rate is about 80 percent in each month.
For each of the responding $n$ units there is a corresponding
random variable, denoted $Y$, that counts the number of responses, 
in the four interviewing attempts. Those with 0 responses are truncated.
Indeed the records for the reason of  0 responses were not accurate, and thus we preferred to ignore/truncate the records that correspond to 0 responses.
We model the distribution of an observed random variable, i.e. conditional 
on $i \in {\bf S}$ by
$$ Y=1+W; \;\; W \sim Binomial(3,p).$$

The above model amounts to assuming that the probability of response
of unit $i$,
is $p_i$ in all of its four investigation attempts, and responses in different months are independent. 
Given a grid
${s}_1,...,{s}_{k}$, for the support of the possible values of $p$, 
a matrix $P=(p_{jk})$
is defined where $p_{jk}=  P(Y=j|p=s_k)=P_{s_k}(1+W=j)$,
$j=1,2,3,4$,
for $W \sim B(3,{s}_k)$. In our analysis we took the grid
0.1, 0.11, 0.12,...,1.  The above induces a matrix $P^*$ in a manner
similar to the previous sections.

Now, $\alpha_{x_l}$, $l=0,1,2$, may be estimated by $\hat{\alpha}^1_{x_l}$ as given  in (\ref{eqn:htalpha})
for a truncation setup.
However, so far we considered only the data from the panel that has four investigations. Indeed the panels that have less investigations
will yield poor estimates of $E( 1/p| X=x_l, I=1)$. Our approach is the following hybrid method. We estimate $E( 1/p| X=x_l, I=1)$,
$l=0,1,2$, based not only on
the data from the current `fourth panel', but, in addition we use 
the data obtained in the four investigation of the three more panels that had their fourth investigation in the previous month,
two months ago, and three months ago, altogether four panels. Let
$m_l$, be the number of items in the {\it currently} investigated four panels, with corresponding $X=x_l$, $x_l=0,1,2$.
Our hybrid approach is to inflate $m_l$, which is based on
the currently investigated four panels, using the estimated $E( 1/p| X=x_l, I=1)$,
$x_l=0,1,2$, which are in turn based on the current as well as `historical' complementary information.
The underlying assumption is that $E( 1/p| X=x_l, I=1)$,  changes slowly in time and thus, estimating it based on a complementary older data,
we still  get at least some bias correction. We proceed by estimating by
$\hat{G}^t$, the joint distribution of $(X,p)$ under truncation.
Finally, we get the  estimator  $$ \hat{\alpha}_{x_{l_0}}=
\frac{ m_{l_0} E_{\hat{G}^t}( 1/p| X=x_{l_0}) }
{ \sum_l  m_l E_{\hat{G}^t}( 1/p| X=x_l)  }.$$

\bigskip

Since the true proportions of the various working statuses
are unknown, we will first demonstrate
the performance of the above estimation method in estimating the following {\it known} true proportions, based on the responses in a given month.

In one case  we  estimate the proportion of males in the population,
which is known to be 0.4853; their proportion in the survey
among responders
is about one percent lower. In the other example we estimate the proportion of the group age 20-39. Their known proportion is 0.397
while their, 
response rate is particularly low, their proportion  among the responders is nearly 3 percent lower than their proportion in the population. 

Each of the following tables 3 and 4 has three lines that correspond
to the  data obtained in Aug/2012, Dec/2012, and April/2013. We took
periods that are four months apart in order not to have overlapping
panels. The general picture persist in other  months.

The columns  True, Naive, and $\hat{\alpha}$, correspond to the true 
population’s proportion,
the sample proportion among responders, and our estimator $\hat{\alpha}$. In each case one may see that $\hat{\alpha}$ corrects the sample proportion
in the right direction.

\bigskip

After gaining some confidence in $\hat{\alpha}$, we will now examine
its estimates in the estimation of the proportion of `Unemployed',
`Employed' and those `Not in Working Force' (NWF).
In the following Table 5 the columns Naive and $\hat{\alpha}$
are as before. The column Bureau gives the estimates of the 
Israel, Central Bureau of Statistics, for the three categories of working statuses. The three parts of the table refer to the three working statuses. The three lines in each part refer to the three months as described before. The Bureau and the $\hat{\alpha}$ estimators `correct'
the naive estimator for Employment and NWF,
in opposite directions  (the official Bureau estimator
involves additional seasonal adjustment that we neglect). The estimator of the bureau is obtained through a method that involves calibration
in a `post-stratification manner'. It seems that the correction of the bureau, of  `Employment' and the `NWF' is in the wrong direction. This is indicated also when imputing missing values based
on their values in months where a response was obtained looking also `into the future'. On the other hand both the Bureau and $\hat{\alpha}$ correct the unemployment naive estimate by increasing it.
This direction of correction of unemployment,
is ,again, supported also by an analysis that involves imputation.   

%\newpage

\begin{table}[ht]
\caption{Comparison of estimates of male's proportion.} % title of Table
\centering % used for centering table
\begin{tabular}{c c c c} % centered columns (4 columns)
\hline\hline %inserts double horizontal lines
   & True & Naive & $\hat{\alpha}$ \\ [0.5ex] % inserts table
%heading
\hline % inserts single horizontal line
{\text Male}&0.4853   & 0.4752 & 0.4822\\
            & 0.4853 & 0.4751 & 0.4819 \\ % inserting body of the table
            &0.4853 &  0.4776 & 0.4842 \\
%{\text UnEmp}& 0.3465 & 0.3576 & 0.3605 \\
% & 31 & 25 & 415 \\
%4 & 35 & 144 & 2356 \\
%{\text NWF} & 0.0454 & 0.0431 & 0.04840\\ [1ex] % [1ex] adds vertical space
\hline %inserts single line
\end{tabular}
\label{table:nonlin2} % is used to refer this table in the text
\end{table}

\begin{table}[ht]
\caption{ Comparison of estimates of  proportion of 20-39
age group.} % title of Table
\centering % used for centering table
\begin{tabular}{c c c c} % centered columns (4 columns)
\hline\hline %inserts double horizontal lines
   & True & Naive & $\hat{\alpha}$  \\ [0.5ex] % inserts table
%heading
\hline % inserts single horizontal line
{\text Age 20-39} & 0.3970 & 0.3664 & 0.3815 \\
                  & 0.3970 & 0.3631 & 0.3984 \\ % inserting body of the table
                 & 0.3970 & 0.3598 &  0.3842\\
%{\text UnEmp}& 0.3465 & 0.3576 & 0.3605 \\
% & 31 & 25 & 415 \\
%4 & 35 & 144 & 2356 \\
%{\text NWF} & 0.0454 & 0.0431 & 0.04840\\ [1ex] % [1ex] adds vertical space
\hline %inserts single line
\end{tabular}
\label{table:nonlin3} % is used to refer this table in the text
\end{table}

\begin{table}[ht]
\caption{Comparison of unemployment estimates.} % title of Table
\centering % used for centering table
\begin{tabular}{c c c c} % centered columns (4 columns)
\hline\hline %inserts double horizontal lines
   & Bureau & Naive & $\hat{\alpha}$ \\ [0.5ex] % inserts table
%heading
\hline % inserts single horizontal line
{\text Emp}
            & 0.6104 & 0.5931 & 0.5761\\
            & 0.6081 & 0.5992 & 0.5910 \\ % inserting body of the table
            &0.6089 & 0.5986 & 0.5881 \\
{\text NWF}&0.3416  & 0.3594 & 0.3748\\
           & 0.3465 & 0.3576 & 0.3605 \\
           &0.3491  & 0.3621 & 0.3720 \\
% & 31 & 25 & 415 \\
%4 & 35 & 144 & 2356 \\
{\text UnEmp} &0.0479  & 0.0475 & 0.0492 \\
              & 0.0454 & 0.0431 & 0.0484\\ 
              &0.0420 & 0.0392 & 0.0399 \\
              %[1ex] % [1ex] adds vertical space
               
\hline %inserts single line
\end{tabular}
\label{table:nonlin1} % is used to refer this table in the text
\end{table}

\newpage

\vspace{3ex}
{\bf \Large References}

\begin{list}{}{\setlength{\itemindent}{-1em}\setlength{\itemsep}{0.5em}}

\item

Benjamini, Y. and Hochberg, Y. (1995). Controlling the false discovery rate: a practical and powerful approach to multiple testing. {\it JRSSB} {\bf 57} No.1, 289-300.

\item
Brown, L.D. and Greenshtein, E. (2009). Non parametric
empirical Bayes and compound decision
approaches to estimation of high dimensional vector of normal
means. {\it Ann. Stat.} {\bf 37}, No. 4, 1685-1704.

\item

Brown L.D., Greenshtein, E. and Ritov, Y. (2013). The Poisson compound decision revisited. {\it JASA.} {\bf 108} 741-749.

%\item

%Efron, B. (2003).  Robbins, Empirical Bayes and Microarrays.
%{\it Ann.Stat.} {\bf 31} No. 2,  366-378.

%\item

%Efron B, and Thisted, R. (1976).  Estimating the number of unseen
%species: How many words did Shakespeare know? {\it Biometrika}
%{\bf 63} 435-447.

\item

Efron, B. (2013). Empirical Bayes modeling, computation and accuracy.
Manuscript.
\item

Greenshtein, E., Park, J., and Ritov, Y. (2008). Estimating the mean of high valued observations in high dimensions. {\it JSTP} {\bf 2} No. 3 407-418.

%\item

%Jiang, W. and Zhang, C.-H. (2009). General maximum likelihood
%empirical Bayes estimation of normal means. {\it Ann. Stat.} {\bf 37}, No 4, 1647-1684.

\item
Koenker, R. and Mizera, I. (2013). Convex optimization, shape constraints, compound decisions and empirical Bayes rules. Manuscript.
\item

Lee, M., Hall, P., Haipeng, S., Marron, J.S., and Tolle, J. (2013).
Deconvolution estimation of mixture distributions with boundaries. (2013). {\it Electronic J. of Stat.} {\bf 7} 323-341.

\item

Lindsay, B. G. (1995). Mixture Models: Theory, Geometry and Applications.
Hayward, CA, IMS.

\item

Lindsay, B. G. and  Roeder, K. (1993).  Uniqueness of estimation and
identifiability in mixture models. {\it Canadian Journal of Stat.}
{\bf 21}, No. 2, 139-147.

%\item

%Link, W. A. (2004). Individual heterogeneity and identifiability in capture recapture models.  {\it Animal biodiversity and conservation} 27.1, 87-91.

\item
Little, R.J.A and Rubin, D.B.  (2002). Statistical Analysis with Missing Data. New York: Wiley.

\item

Karatzoglou, A., Smola, A., Hornik, K., and Zeleis, A., (2004).
An S4 package for kernel methods in R. {\it Journal of Statistical Software}
{\bf 11}, No. 9, 1-20.

\item

Kiefer and Wolfowitz (1956). Consistency of the maximum likelihood estimator in the presence  of infinitely many incidental parameters.
{\it Ann.Math.Stat.} {\bf 27} No. 4, 887-906.

%\item

%Norris, J. L. and Pollock, K. H. (1996). Nonparametric MLE under two closed capture-recapture
%models with heterogeneity. {\it Biometrics }, {\bf 52}, 639�649.

\item
Sharon L. Lohr (1999). \emph{Sampling Design and Analysis}. Brooks/Cole publishing company.

%\item

%Vapnik, N. V. (1998). Statistical Learning Theory. Wiley, New York.

\item

Zhang, C-H. (2005). Estimation of sums of random variables: Examples and information bounds. {\it Ann. Stat.} {\bf 33} No.5. 2022-2041.

\end{list}

\end{document}